\documentclass[12pt]{amsart}
\usepackage{amsmath, amsthm, euscript}
\swapnumbers \theoremstyle{plain}
\newtheorem{thm}{Theorem}[section]

\newtheorem{cor}[thm]{Corollary}
\numberwithin{equation}{section}

\theoremstyle{definition}
\newtheorem{setup}[thm]{Setup}
\newtheorem{remark}[thm]{Remark}
\newtheorem{K_n}[thm]{Quantized $\k$-algebras $K_n$}
\newtheorem{quan-matr}[thm]{Quantum Matrices}
\newtheorem*{Ack}{Acknowledgment}

\theoremstyle{remark}
\newtheorem{note}[thm]{}

\def\rKdim{\operatorname{rKdim}}
\def\gr{\operatorname{gr}}
\def\rgl{\operatorname{rgl}}

\def\m{\mathfrak{m}}
\def\t{\tau}
\def\d{\delta}
\def\s{\sigma}

\def\tl{\t_l}
\def\dl{\d_l}
\def\y{Y_{l-1}}

\def\tlm{\t_{l m}}
\def\dlm{\d_{l m}}

\def\k{\mathbf{k}}
\def\N{\mathbb{N}}

\def\linei{\underline{i}}
\def\linej{\underline{j}}

\begin{document}
\title{Completions of Quantum Coordinate Rings}

\author{Linhong Wang}

\address{Department of Mathematics\\
        Temple University\\
        Philadelphia, PA 19122-6094}

\email{lhwang@temple.edu}


\begin{abstract}
Given an iterated skew polynomial ring $C[y_1;\t_1,\d_1]\ldots [y_n;\t_n,\d_n]$ over a complete local ring $C$ with maximal ideal $\m$, we prove, under suitable assumptions, that the completion at the ideal $\m + \left\langle y_1,y_2,\ldots,y_n\right\rangle$ is an iterated skew power series ring. Under further conditions, the completion becomes a local, noetherian, Auslander regular domain. Applicable examples include quantum matrices, quantum symplectic spaces, and quantum Euclidean spaces.
\end{abstract}

\maketitle
\section{Introduction}
Let $R$ be a ring equipped with a skew derivation $(\t,\d)$. The skew power series ring $R[[y;\tau]]$, when $\delta = 0$, is a well known, classical object (cf. \cite{Cohn}, \cite{McCon-Robs}). The skew power series ring $R[[y;\t,\d]]$, when $\delta \neq 0$, has more recently appeared in quantum algebras (cf. \cite [\S 4] {Kash}, \cite [\S 4]{Koo}) and in noncommutative Iwasawa theory (cf. \cite{Schn-Ven}, \cite{Ven}). In this paper, we study iterated skew power
series rings as completions of iterated skew polynomial rings. Our approach builds on the work of Venjakob in \cite{Ven}.

Our main result can be stated as follows: \emph{Let $$R_n=C[y_1;\t_1,\d_1]\ldots [y_l;\t_l,\d_l] \ldots [y_n;\t_n,\d_n]\quad (n\geq 1)$$ be an iterated skew
polynomial ring, where $C$ is a complete local ring with maximal ideal $\m$, and where $C$ is stable under each skew derivation $(\tl,\dl)$. For
each $1\leq l \leq n$, set $I_{l-1}=\m + \left\langle y_1,\ldots, y_{l-1} \right\rangle$, and assume that $\tl(I_{l-1}) \subseteq I_{l-1}$,
$\dl(R_{l-1})\subseteq I_{l-1}$, and $\dl(I_{l-1})\subseteq I_{l-1}^2$. Then there exists an iterated skew power series ring
\[S_n=C[[y_1;\hat{\t}_1,\hat{\d}_1]] \ldots [[y_l;\hat{\t}_l,\hat{\d}_l]] \ldots [[y_n;\hat{\t}_n,\hat{\d}_n]],\]
such that $\hat{\t}_l \mid_{R_{l-1}} =\tl$ and $\hat{\d}_l \mid_{R_{l-1}} =\dl$, for $1\leq l \leq n$. Moreover, $S_n$ is the completion of $R_n$ at
the ideal $\m + \left\langle y_1,\ldots, y_n \right\rangle$.}

The paper is organized as follows: Section 2 reviews some preliminary results and proves the main result. Section 3 applies the main result to certain quantum coordinate rings, including quantum matrices, quantum symplectic spaces, and quantum Euclidean spaces.

Throughout, all rings are unital.

\section{Main Result}
Let $R$ be a ring, $\t$ a ring endomorphism of $R$ and $\d$ a left $\t$-derivation, that is, $\d:\, R\rightarrow R$ is an additive map for which
$\d(rs)=\t(r)\d(s)+\d(r)s$ for all $r,\,s \in R$. The pair of maps $(\t,\, \d)$ is called a \emph{skew derivation} on $R$. To start, we recall the structure of the skew
power series ring in one variable, following Venjakob \cite{Ven}.

\begin{note}\label{2r1}
Let $S$ be the additive group of formal power series in $y$, \[\sum_ir_iy^i= \sum_{i=0}^\infty r_i y^i,\] with coefficients $r_i$ in $R$. Using the commutation rule $yr=\t(r)y+\d(r)$, for $r\in R$, we wish to write the product of two arbitrary elements in $S$ as
\[ \left(\sum_ir_iy^i\right)\left(\sum_js_jy^j\right)=\sum_n^{\infty}
\sum_{j=0}^n \sum_{i=n-j}^{\infty}r_i(y^is_j)_{n-j}y^n, \]
where each $(y^n r)_i$, for $0\leq i\leq n$, denotes an element in $R$ such that \[y^n r=\sum_{i=0}^n(y^n r)_i y^i,\] for $n\geq 0$. However, it is
not always the case that  \[\sum_{j=0}^n \sum_{i=n-j}^{\infty}r_i(y^is_j)_{n-j}\] is well defined in $R$. If, under some additional restrictions
(see \ref{2r3}), the multiplication formula is well defined for any two power series in $S$, we will say that $S$ is a \emph{well-defined skew power
series ring}, and write $S=R[[y;\t,\d]]$.
\end{note}

\begin{note}\label{2r2}
By a \emph{local ring} we will always mean a ring $R$ such that the quotient ring by the Jacobson radical $J(R)$ is simple artinian. In particular,
a local ring has a unique maximal ideal which is equal to the Jacobson radical. Let $R$ be a local ring with maximal ideal $\m$.  We will always equip $R$ with the $\m$-adic topology. By the associated graded ring $\gr R$, we
will always mean with respect to the $\m$-adic filtration, that is: \[\gr R= R/\m \oplus \m/\m^2 \oplus \cdots\] We will refer to  $R$ as a
\emph{complete local ring} if $R$ is also \emph{complete} (i.e., Cauchy sequences converge in the $\m$-adic topology) and \emph{separated} (i.e.,
the $\m$-adic topology is Hausdorff).
\end{note}

\begin{note} \label{2r3}
Let $R$ be a complete local ring with maximal ideal $\m$ and with skew derivation $(\t,\d)$. As in \cite{Ven}, we assume that
$\t(\m)\subseteq\m,\; \d(R)\subseteq\m\;\text{and}\; \d(\m)\subseteq \m^2$. In \cite[Lemma 2.1]{Ven}, Venjakob proved, under these assumptions, that $S=R[[y;\t,\d]]$ is
a well-defined skew power series ring. The following properties of $S$ are also proved in, or easily deduced from, Venjakob's work in \cite[section 2]{Ven}.\\
(i) Any element $\sum_i r_i y^i$ is a unit (in $S$) if and only if the constant term $r_0$ is a unit in $R$. In particular, any element in
$1-\left\langle \m,\; y\right\rangle$ is a unit, and so the Jacobson radical $J(S)=\left\langle \m,\; y\right\rangle$. Hence, in view of the
isomorphism $S/J(S) \cong R/\m$, $S$ is a local ring.\\
(ii) The $\left\langle \m,\; y\right\rangle$-adic filtration on $S$ is complete and separated.\\
(iii) There is a canonical isomorphism $\gr S \cong (\gr R)[\bar{y};\,\bar{\t}]$. Assume further that $\bar{\t}$ is an automorphism. Then, $S$ is right (respectively left) noetherian if $\gr R$ is right (respectively left) noetherian, $S$ is a domain if $\gr R$ is a domain, and $S$ is Auslander regular if the same holds for $\gr R$; see
\cite[Corollary 2.10]{Ven} (cf. \cite[Chap. III, Theorem 2.2.5]{Li-Oystaeyen}, \cite[Chap. III, Theorem 3.4.6 (1)]{Li-Oystaeyen}).\\
(iv) Now suppose that $\gr R$ is right noetherian and that $\bar{\t}$ is an automorphism. Concerning right global dimension, it holds that $\rgl S \leq \rgl\, \gr R +1$. As far as right Krull dimension is concerned, $\rKdim \gr(S)$ $= \rKdim \gr R +1$, by \cite[15.19]{Goodearl-Warfield}. Moreover, $\rKdim S \leq \rKdim\gr S$, as $S$
is a complete filtered ring and $\gr S$ is right noetherian; see \cite[D.IV.5]{NasVOy}. Therefore $\rKdim S \leq \rKdim \gr R +1$.
\end{note}

\begin{note}\label{2r4}
Contained within $S$ is the skew polynomial ring $T=R[y;\,\t,\d]$. Following \ref{2r3} (ii), both $S$ and $T$ are endowed with a Hausdorff
$\left\langle \m,\,y\right\rangle$-adic topology. Of course, $T$ is a dense subring of $S$ in this topology. Therefore, $S$ is the completion of $T$
with respect to the $\left\langle \m,\,y\right\rangle$-adic filtration, following \cite[3.3.5]{AGM-topo}.
\end{note}

The remainder of this section is devoted to the main result. First we set up a suitable iterated skew polynomial ring. Then we construct an iterated skew power series ring, by extending skew derivations.

\begin{setup}\label{2setup}
Let $C$ be a complete local ring with maximal ideal $\m$. Set $R_0=C$, and let
\[R_n=C[y_1;\t_1,\d_1] \ldots [y_l;\tl,\dl] \ldots
[y_n;\t_n,\d_n]\] be an iterated skew polynomial ring with skew derivations $(\tl,\dl)$ of $R_{l-1}$, for $1\leq l \leq n$. For each $1\leq l \leq n$, set
\[ I_{l-1}=\m + \left\langle y_1,\ldots, y_{l-1} \right\rangle \subseteq R_{l-1},\]
and assume that
\[ \tl(I_{l-1}) \subseteq I_{l-1},\quad \dl(R_{l-1})\subseteq I_{l-1}, \quad\text{and}\quad \dl(I_{l-1})\subseteq I_{l-1}^2.\]
\end{setup}

We will also need the following notations.

\begin{note}\label{2-notation}
(i) Let $1\leq l \leq n+1$. A nonzero monomial $c_{i_1,\ldots,i_{l-1}} y_1^{i_1}\cdots y_{l-1}^{i_{l-1}}$ in $R_{l-1}$ is said to be in \emph{normal form}.
We will write $$c_{\linei} \y^{\linei}$$ for  $c_{i_1,\ldots,i_{l-1}} y_1^{i_1}\cdots y_{l-1}^{i_{l-1}}$, where  $\linei= (i_1,\ldots,i_{l-1}) \in
\N^{l-1}$.\\
(ii) We now introduce the notion of degree we will use for nonzero monomials in normal form. Let $1\leq l \leq n+1$, and let $c_{\linei} \y^{\linei}\in R_{l-1}$. Then there exists an
integer $k$ largest such that $c_{\linei}\in \m^{k}$. Set
$$s(c_{\linei}, \linei)=k + i_1+i_2+\cdots+i_{l-1}.$$
We will refer to  $s(c_{\linei}, \linei)$ as the \emph{degree} of $c_{\linei} \y^{\linei}$.\\
(iii) Let $1\leq l \leq n$, and let $c_{\linei} \y^{\linei}$ and $d_{\linej} \y^{\linej}$ be two nonzero monomials in $R_{l-1}$. Then $c_{\linei} \y^{\linei} \cdot d_{\linej}
\y^{\linej}$ is $0$ or a sum of monomials each with degree $\geq s(c_{\linei}, \linei) + s(d_{\linej}, \linej)$. An inductive argument shows that each of the polynomials $\tl\left(c_{\linei} \y^{\linei}\right)$ and $\dl\left(c_{\linei} \y^{\linei}\right)$ is $0$ or a finite sum of monomials each with degree $\geq s(c_{\linei}, \linei)$. \\
(iv) Let $1\leq l \leq n$. By a \emph{formal power series} in $y_1,\ldots, y_l$ over $C$, we will mean an infinite series \[f=\sum_{\linei} c_{\linei} Y_l^{\linei},\]
where the $c_{\linei}$ are elements in $C$ and where $\linei \in \N^l$. Note that each monomial $c_{\linei} Y_l^{\linei}$ is in normal form. The set
of all formal power series in $y_1,\ldots, y_l$ over $C$ forms an abelian group, which we will denote as $A_l$.
\end{note}
\begin{note}\label{2-wdefined} Let $1\leq l \leq n$.\\
(i) Given a power series $f=\sum_{\linei} c_{\linei} Y_l^{\linei} \in A_l$, we can always write
\[f= \sum_{k=0}^\infty \sum_{s(c_{\linei}, \linei)=k} c_{\linei} Y_l^{\linei},\]
after regrouping the monomials appearing in $f$ (if necessary). Note that for each $k$ the sum \[\sum_{s(c_{\linei}, \linei)=k} c_{\linei}
Y_l^{\linei}\] is finite and possibly equal to 0.\\
(ii) On the other hand, let
\[g=G_0 + G_1 + \ldots+ G_k+ \ldots ,\]
where each $G_k$ is $0$ or a finite sum of monomials in $R_l$ all with degree $k$. Then $g$ is a well-defined (in the above sense)
formal power series in $A_l$. To see this, suppose that
\[G_k=\sum_{\linej \in M_k } c_{\linej}^{(k)} Y_l^{\linej},\]
where $c_{\linej}^{(k)}\in C$ and where $M_k \subseteq \N^l$, for k=0,1,\ldots. We will set $c_{\linei}^{(k)}=0$ when $\linei \notin M_k$. Now, for
a fixed $\linej$, the sum
$$ c_{\linej}^{(0)} + c_{\linej}^{(1)} +\ldots + c_{\linej}^{(k)}+ \ldots $$
might contain infinitely many terms. But each $c_{\linej}^{(k)}$ is such that the degree of $c_{\linej}^{(k)} Y_l^{\linej}$ is equal to
$k$. Hence, the preceding sum  is convergent in the $\m$-adic topology. Therefore,
$$g=G_0 + G_1 + \ldots+ G_k+ \ldots
=\sum_{\linej \in \cup M_k} \left( c_{\linej}^{(0)} + c_{\linej}^{(1)} +\ldots + c_{\linej}^{(k)}+ \ldots\right)Y_l^{\linej}$$
is a formal power series in $A_l$ with all coefficients in $C$ well defined.
\end{note}

\begin{thm}\label{2t}
Retain notations and assumptions in \ref{2setup}. Let $S_0=C$. Then there exists an iterated skew power series ring
\[S_n=C[[y_1;\hat{\t}_1,\hat{\d}_1]] \ldots [[y_l;\hat{\t}_l,\hat{\d}_l]] \ldots [[y_n;\hat{\t}_n,\hat{\d}_n]],\]
where each $(\hat{\t}_l,\hat{\d}_l)$ is a skew derivation on $S_{l-1}$ with $\hat{\t}_l \mid_{R_{l-1}} =\tl$ and $\hat{\d}_l \mid_{R_{l-1}} =\dl$,
for $1\leq l \leq n$. Moreover, $S_n$ is a complete local ring with maximal ideal $\m_n=\m + \left\langle y_1,\ldots, y_n \right\rangle$. (We will
refer to $S_n$ as the \emph{power series extension} of $R_n$.)
\end{thm}

\begin{proof}
Following \ref{2r3}, the ring $C[[y_1;\t_1,\d_1]]$ is well defined and we may take $S_1=C[[y_1;\t_1,\d_1]]$. In the notation of \ref{2-notation},
$S_1$ is the abelian group $A_1$ equipped with a well-defined multiplication restricting to the original multiplication in $R_1$. Our goal is to
show that each abelian group $A_l$ becomes an iterated skew power series ring. In the first step of the proof, we extend the pair of maps $\t_l$ and
$\d_l$ to $A_{l-1}$ for all $ 1 < l \leq n$. Then, by induction, we will show that each $(\t_l,\d_l)$ extends to a skew derivation on $S_{l-1}$ and
that each $A_{l}$ forms a ring $S_l$ .\smallskip

To start, let $f=\sum_{\linei} c_{\linei} \y^{\linei}$ be a power series in $A_{l-1}$. As in \ref{2-wdefined} (i), we can write
\[f=\sum_{k=0}^{\infty} F_k,\quad \quad \text{where} \quad F_k:=  \sum_{s(c_{\linei}, \linei) =k} c_{\linei} \y^{\linei}\; \,\text{(possibly equal to 0)}\]
Our goal now is to extend $\t_l$ and $\d_l$ to $A_{l-1}$. For $k=0,1,2,\ldots$, if $\t_l(F_k) \neq 0$, then \[\t_l(F_k)=\sum_{\linej \in T_k} t_{\linej}^{(k)} \y^{\linej}\] for some
subset $T_k \subseteq \N^{l-1}$ and some $t_{\linej}^{(k)}\in C$.
Next, let
\[G_m= \sum_{k=0}^\infty \sum_{\linej\in N_{m,k}} t_{\linej}^{(k)} \y^{\linej},\] where
\[N_{m,k}=\{\linej \in T_k \mid \text{the degree of } t_{\linej}^{(k)} \y^{\linej}\; \text{is}\; m\}.\] In other words, we regroup the monomials appearing in $\sum_k \t_l(F_k)$ by their degrees. Then
\begin{equation*}\label{tF-G}
\t_l(F_0) +\t_l(F_1)+\ldots + \t_l(F_k)+\ldots\; =\; G_0+G_1+\ldots +G_m+\ldots
\end{equation*}
It follows from \ref{2-notation} (iii) that any nonzero $\t_l(F_k)$ is a finite sum and that each $t_{\linej}^{(k)} \y^{\linej}$ has degree $\geq k$. Hence each $G_m$ is a finite sum by the construction. Recall from \ref{2-wdefined} (ii) that $$G_0+G_1+\ldots +G_m+\ldots$$ is a formal power series in $A_{l-1}$. Therefore,
$$\sum_{k=0}^{\infty} \t_l \left( F_k \right)\in A_{l-1}.$$ Using the same argument (replacing $\t_l$ with $\d_l$), we also have
$$\sum_{k=0}^{\infty} \d_l \left( F_k\right) \in A_{l-1}.$$
Then, for $1\leq l \leq n$ and $f= \sum_{\linei} c_{\linei} \y^{\linei} \in A_{l-1}$, we extend $\tl$ and $\dl$ by setting the following maps
\begin{equation}\label{t-d}
\hat{\t}_l (f)= \sum_{k=0}^{\infty} \t_l \left( \sum_{s(c_{\linei}, \linei) =k} c_{\linei} \y^{\linei} \right)\quad \text{and}\quad
\hat{\d}_l (f)= \sum_{k=0}^{\infty} \d_l \left( \sum_{s(c_{\linei}, \linei) =k} c_{\linei} \y^{\linei} \right).
\end{equation}
It is clear that $\hat{\t}_l \mid_{R_{l-1}} =\tl \quad \text{and} \quad \hat{\d}_l
\mid_{R_{l-1}} =\dl$, for all $1\leq l \leq n$. \medskip

Now, let $n \geq 2$. Assume that the abelian group $A_{n-1}$ is a well-defined power series ring, which we will denote as $S_{n-1}$, and also assume that
$S_{n-1}$ is a complete local ring with maximal ideal $\m_{n-1}=\m + \left\langle y_1,\ldots, y_{n-1} \right\rangle$. Next we show that
$(\hat{\t}_n,\,\hat{\d}_n)$, from (\ref{t-d}), is a skew derivation on $S_{n-1}$; that is, $\hat{\t}_n$ is an automorphism of $S_{n-1}$ and
$\hat{\d}_n$ is a left $\hat{\t}_n$-derivation.\medskip

Let $t$ be a positive integer. Choose two arbitrary elements $a$ and $b$ in $S_{n-1}$. Write $a=a_t + a_t'$ and $b=b_t + b_t'$, where $a_t$
(respectively $b_t$) is the sum of the monomials appearing in $a$ (respectively $b$) with degree $\leq t$. Then it follows from (\ref{t-d}) that
\[ \hat{\t}_n (a)= \hat{\t}_n (a_t)+ \hat{\t}_n (a_t') \quad \text{and} \quad \hat{\t}_n (b)= \hat{\t}_n (b_t)+\hat{\t}_n (b_t').\]
Therefore, we have
\[ \hat{\t}_n (ab)= \t_n \left(a_t \cdot b_t \right) + \hat{\t}_n (a_t' \cdot b_t+a_t \cdot b_t' +a_t '\cdot b_t'), \quad \text{and}\]
\[ \hat{\t}_n (a) \cdot \hat{\t}_n (b)= \t_n(a_t) \cdot \t_n(b_t) + \hat{\t}_n (a_t') \cdot \hat{\t}_n (b_t)+ \hat{\t}_n (a_t) \cdot
\hat{\t}_n (b_t') + \hat{\t}_n (a_t') \cdot \hat{\t}_n (b_t').\]
Note that $\t_n(a_t \cdot b_t)=\t_n(a_t) \cdot \t_n(b_t)$. It follows from \ref{2-notation} (iii) that
\[\hat{\t}_n (ab)-\hat{\t}_n (a) \cdot \hat{\t}_n (b) \in \m_{n-1}^{t+1}.\]
Let $t \rightarrow \infty$, then it follows from the completeness of $S_{n-1}$ that
\[\hat{\t}_n (ab)=\hat{\t}_n (a) \cdot \hat{\t}_n (b).\]
Using the same argument (replacing $\hat{\t}_n$ with $\hat{\d}_n$), we can get
 \[\hat{\d}_n (ab)= \hat{\d}_n (a) b + \hat{\t_n}(a) \hat{\d}_n(b).\] Therefore $(\hat{\t}_n,\,\hat{\d}_n)$
is a skew derivation on $S_{n-1}$.\medskip

In view of the assumptions in \ref{2setup} and (\ref{t-d}), we see that
\[\hat{\t}_n(\m_{n-1}) \subseteq \m_{n-1}, \quad \hat{\d}_n(S_{n-1}) \subseteq \m_{n-1},\quad \text{and}\quad \hat{\d}_n(\m_{n-1}) \subseteq
\m_{n-1}^2.\] Following \ref{2r3} (i) (ii), the skew power series ring $S_n=S_{n-1}[[y_n;\,\t_n,\,\d_n]]$ is well defined, and $S_n$ is a complete local ring
with maximal ideal $\m_n=\m + \left\langle y_1,\ldots, y_{n} \right\rangle$. This completes the inductive step. The theorem is proved by induction.
\end{proof}

The following is a consequence of \ref{2r3}, \ref{2r4} and Theorem \ref{2t}.

\begin{cor}\label{2c1}
{\rm (i)} The power series extension $S_n$ in \ref{2t} is the completion of $R_n$ with respect to the ideal $\m_{n}=\m + \left\langle y_1,\ldots, y_{n} \right\rangle$. Any power series in $S_n$ is a unit (in $S_n$) if and only if its constant term is a unit in $C$.\\
{\rm (ii)} The associated graded ring $\gr S_n$ is isomorphic to an iterated skew polynomial ring $(\gr C)[y_1;\,\bar{\t}_1] \, \ldots \,[y_n;\,\bar{\t}_n]$.\\
{\rm (iii)} Assume further that $\bar{\t}_1,\ldots,\bar{\t}_n$ are automorphisms. If $\gr C$ is a domain, $S_n$ is a domain. If $\gr C$ is right (respectively left) noetherian, so is $S_n$. If $\gr C$ is Auslander regular, then $S_n$ is also Auslander regular.\\
{\rm (iv)} Suppose that $\gr C$ is right noetherian and that $\bar{\t}_1,\ldots,\bar{\t}_n$ are automorphisms. Then it holds that $\rKdim S_n \leq \rKdim \gr C +n$ and $\rgl S_n \leq \rgl \gr C +n$.
\end{cor}

\section{Examples}
Throughout, let $\k$ be a field.
\begin{quan-matr}\label{3e1}
Let $\mathcal{O}_{\lambda,\,\mathbf{p}}(M_n(\k))$ be the multiparameter quantum coordinate ring of $n\times n$ matrices over $\k$, as studied in
\cite{Artin-S-Tate} (cf. e.g., \cite{Brown-Goodearl}). Here $\mathbf{p}=(p_{i j})$ is a multiplicatively antisymmetric $n\times n$ matrix over $\k,$ and $\lambda$ is a nonzero
element of $\k$ not equal to $1$. Further information of this algebra can be found in \cite{Brown-Goodearl}. As shown in \cite{Artin-S-Tate},
$\mathcal{O}_{\lambda,\,\mathbf{p}}(M_n(\k))$ can be presented as a skew polynomial ring
\[\k[y_{11}]\, [y_{1\,2};\t_{1\,2}]\, \cdots\,
[y_{lm};\,\tlm,\dlm]\, \cdots \, [y_{nn};\t_{nn},\d_{nn}].\]
Each $(\tlm,\dlm)$ is a skew derivation as follows:
\begin{align*}
\tlm(y_{ij}) &=\left\{
  \begin{array}{ll}
    p_{li}p_{jm}y_{ij}, & \mathrm{when}\;l\geq i \;\mathrm{and}\; m>j, \\
    \lambda p_{li}p_{jm}y_{ij}, & \mathrm{when}\;l>i \;\mathrm{and}\; m\leq j,\\
    \end{array}
\right. \\
\dlm(y_{ij}) &=\left\{
  \begin{array}{ll}
    (\lambda -1)p_{li}y_{im}y_{lj}, & \mathrm{when}\;l > i\;
    \mathrm{and}\; m>j, \\
    0, & \mathrm{otherwise}.\\
    \end{array}
\right.
\end{align*}
It is not hard to see that these skew derivations satisfy the assumptions in \ref{2setup}. Hence, by Theorem \ref{2t}, the power series extension
of $\mathcal{O}_{\lambda,\,\mathbf{p}}(M_n(\k))$ is the iterated skew power series ring
\[\k[[y_{11}]]\, [[y_{1\,2};\hat{\t}_{1\,2}]]\, \cdots\,
[[y_{lm};\,\hat{\t}_{lm},\hat{\d}_{lm}]]\, \cdots \, [[y_{nn};\hat{\t}_{nn},\hat{\d}_{nn}]],\] where each extended skew derivation is defined as in
(\ref{t-d}). Also note that each $\t_{lm}$ acts by nonzero scalar multiplication on the generators, and so each $\bar{\t}_{lm}$ is an automorphism. It now follows from \ref{2c1} that the preceding power series completion is a local, noetherian, Auslander regular domain.

\end{quan-matr}

\begin{K_n} \label{3e2}
There are other well-known quantum coordinate rings, for example coordinate rings of quantum symplectic space and quantum Euclidean $2n$-space (see, e.g., \cite{Brown-Goodearl}). Horton introduced a class of algebras, denoted $K_{n,\,\Gamma}^{P,\,Q}(\k)$ or more briefly $K_n$, that includes coordinate rings of both
quantum symplectic space and quantum Euclidean $2n$-space; see \cite{Horton}. To describe this class of algebras, let $P,\,Q \in (\k^{\times})^n$ such
that $P=(p_1,\ldots,\,p_n)$ and $Q=(q_1,\,\ldots,\,q_n)$ where $p_i\neq q_i$ for each $i\in \left\lbrace 1,\,\ldots,\,n\right\rbrace$. Further, let
$\Gamma=(\gamma_{i,j})\in M_n(\k^{\times})$ with $\gamma_{j,i}=\gamma_{i,j}^{-1}$ and $\gamma_{i,i}=1$ for all $i,\,j$. Then, as in \cite{Horton},
$K_{n,\,\Gamma}^{P,\,Q}(\k)$ is generated by $x_1,\, y_1,\,\ldots, x_n,\, y_n$ satisfying certain relations determined by $P,\,Q$ and $\Gamma$. This
algebra can be presented as an iterated skew polynomial ring,
\[
\k[x_1][y_1;\,\t_1][x_2;\,\s_2][y_2;\,\t_2,\,\d_2] \cdots[x_n;\,\s_n][y_n;\,\t_n,\,\d_n];\]
see \cite[Proposition 3.5]{Horton}. Automorphisms $\s_i,\;\t_i$ and $\t_i$-derivations $\d_i$ are defined as follows:
\[
\begin{array}{ll}
\s_i (x_j) = q_j^{-1} p_i \gamma_{i,j} x_j  & 1 \leq j \leq i-1,\\
\s_i (y_j) = q_j \gamma_{j,i} y_j      & 1 \leq j \leq i-1,\\
\t_i (x_j) = p_i^{-1} \gamma_{j,i} x_j      & 1 \leq j \leq i-1,\\
\t_i (y_j) = \gamma_{i,j} y_j               & 1 \leq j \leq i-1,\\
\t_i (x_i) = q_i^{-1} x_i,\\
\d_i(x_j)=0                                 & 1 \leq j \leq i-1,\\
\d_i(y_j)=0                                 & 1 \leq j \leq i-1,\\
\d_i(x_i)=-q_i^{-1} \sum_{l<i} (q_l - p_l)y_l x_l.
\end{array}
\]
Note that these automorphisms and derivations give quadratic relations, and so, by Theorem \ref{2t}, $K_n$ has the power series extension
\[
\k[[x_1]][[y_1;\,\hat{\t}_1]][[x_2;\,\hat{\s}_2]]
[[y_2;\,\hat{\t}_2,\,\hat{\d}_2]]\cdots
[[x_l;\,\hat{\s}_l]][[y_l;\,\hat{\t}_l,\,\hat{\d}_l]]
\cdots[[x_n;\,\hat{\s}_n]][[y_n;\,\hat{\t}_n,\,\hat{\d}_n]],\] where the extended skew
derivations are defined as in (\ref{t-d}). Again, it follows from Corollary \ref{2c1} that this completion is a local, noetherian, Auslander regular domain.
\end{K_n}

Moreover, comparing with \ref{2c1} (iv), the dimensions of the power series completions in \ref{3e2} can be more precisely determined as follows:

\begin{note}
Let $E$ be an algebra in the class $K_n$, and $\hat{E}$ be the power series completion of $E$ with respect to the ideal $\left\langle x_1,y_1,\ldots, x_n,y_n \right\rangle$. From \ref{3e2}, we see that, among the defining commutation relations, nonzero derivations only occur in the following cases:
\[y_ix_i=\t_i(x_i)y_i+\d_i(x_i),\quad \text{for} \; \, i=2,\ldots, n.\]
Also note that $\d_i(x_i)\in I_{i-1}=\langle x_1,y_1,\ldots , x_{i-1},y_{i-1}\rangle$. Hence, the set of generators $\lbrace x_1,y_1,\ldots, x_n,y_n\rbrace$ forms a \emph{regular normalizing set} (see \cite[Definition 1.1]{Walker}). Since $J(\hat{E})=\left\langle x_1,y_1,\ldots, x_n,y_n \right\rangle$, it now follows from \cite[Theorem 2.7]{Walker} that the Krull dimension, classical Krull dimension and global dimension of $\hat{E}$ are all equal to $2n$.
\end{note}

\begin{remark} \label{nilpotent}
For the quantum coordinate rings and quantum algebras in \ref{3e1} and \ref{3e2}, it is well known that the derivations $\d_{lm}$ and $\d_l$ are locally
nilpotent.  In \cite{Cauchon}, using this fact (and other assumptions), Cauchon constructed the ``Derivation-Elimination Algorithm". But, for power series completions of these examples, the extended derivations $\hat{\d}_{lm}$ and $\hat{\d}_l$ are not locally nilpotent.
\end{remark}

\begin{Ack}
This work forms a portion of the author's Ph.D. thesis at Temple University. She is grateful to her thesis advisor E. S. Letzter for his invaluable guidance and generous help. The author is also thankful to the referee for careful reading and helpful comments.
\end{Ack}

\end{document}